\documentclass[journal]{IEEEtran}
\usepackage{amsmath}
\ifCLASSINFOpdf  
\usepackage[pdftex]{graphicx}
\else
\usepackage[dvips]{graphicx}
\fi
\usepackage{epstopdf}
\usepackage[ruled]{algorithm2e} 
\usepackage{booktabs}
\usepackage{multirow}
\usepackage{rotating}
\usepackage{color}
\usepackage{ulem}
\usepackage{color}
\usepackage{cite}
\allowdisplaybreaks[4]
\ifCLASSINFOpdf
\else
\fi
\begin{document}
\bibliographystyle{plain}
%
\title{Fast Green Function Evaluation for \\Method of Moment}
%
%
%

\author{Shunchuan~Yang,~\IEEEmembership{Member,~IEEE,}
        Donglin~Su,~\IEEEmembership{Member,~IEEE}
\thanks{Manuscript received xxx; revised xxx.}
\thanks{This work was supported in part by the National Natural Science Foundation of China through Grant 61801010, Grant 61727802, the ``$13$th Five-Year" Equipment Pre-research Fund through Grant 61402090602 and Beijing Natural Science Foundation through Grant 4194082.}
\thanks{All authors are with the School of Electronic and Information Engineering, Beihang University, Beijing, 10083, China (e-mail: scyang@buaa.edu.cn, sdl@buaa.edu.cn).}
}

%
%

\markboth{ }%
{Shell \MakeLowercase{\textit{et al.}}: Bare Demo of IEEEtran.cls for IEEE Journals}
%



\maketitle

\begin{abstract}
In this letter, an approach to accelerate the matrix filling in method of moment (MOM) is presented. Based on the fact that the Green function is dependent on the Euclidean distance between the source and the observation points, we constructed an efficient adaptive one-dimensional interpolation approach to fast calculate the $Exp$ type function values. In the proposed method, several adaptive interpolation tables are constructed based on the maximum and minimum distance between any two integration points with local refinement near zero function values to minimize the relative error. An efficient approach to obtain the sampling points used in the interpolation phase is carefully designed. Then, any function values can be efficiently calculated through a linear interpolation method for $Exp$ and a Lagrange polynomial interpolation method for the Green function. In addition, the error bound of the proposed method is rigorously investigated. The proposed method can be quite easily integrated into the available MOM codes for different integration equation (IE) formulations with few efforts. Comprehensive numerical experiments validate its accuracy and efficiency through several IE formulations. Results show that over 20\% efficiency improvement can be achieved without sacrificing the accuracy.
\end{abstract}


%
\IEEEpeerreviewmaketitle

\section{Introduction}
%
%
%
%
\IEEEPARstart{T}{He} method of moment (MOM) is widely used to extract composite parameters of large-scale integrated circuits (ICs) \cite{FastImp}, solve electrically large, multiscale scattering problems \cite{MLFMP} due to its unknowns residing on the interface of different homogenous media. Therefore, the count of unknowns is much smaller than that of other differential equation based methods, like the finite-difference time-domain (FDTD) method \cite{FDTD}, the finite element method (FEM) \cite{FEM}, which require volumetric discretization.

However, investigations found that matrix filling in MOM is quite time-consuming and is one of the bottlenecks to prohibit MOM from solving practical problems, which is directly related to calculate
\begin{equation}{\label{EXPRR}}
Exp(-jkr), Exp(-jkr)/r.
\end{equation}
When constructing the impedance matrix, one can directly evaluate those analytical expressions. However, it's reported that it is time-consuming to calculate $Exp$ type function, especially for complex numbers \cite{EXP, EXP2, EXP3}, which is the exact scenario in MOM for the time-harmonic electromagnetic simulations.

In this paper, we show that direct evaluation of those analytical formulations is not optimal since a large number of time-consuming and redundant evaluation of $Exp$ function may exist. We proposed an efficient adaptive interpolation approach to fast calculate  (\ref{EXPRR}). Through careful design of interpolation tables and an fast approach to retrieve the sampling points, we can efficiently calculate their function values with controllable accuracy and significant efficiency improvement.

This paper is organized as follows. In Section II, some key observations related to the integral identities used in MOM are presented and the proposed approach is illustrated in detail. In Section III, its error bound is rigorously analyzed. In Section IV, its accuracy and efficiency are comprehensively investigated through several experiments. At last, we draw some conclusions in Section V.

\section{Method}
\subsection{Observations of Integral Identities Used in MOM}
In this paper, the singular cancellation technique along with the angular transformation \cite{GIBC} is used to demostrate the proposed method. However, it can be easily extended to other integration techniques without any difficulties. Four integrals, Equ. (56) in \cite{GIBC}, are required to be efficiently evaluated through  (\ref{EXPRR}). To efficiently construct the impedance matrix in MOM, several types of integrand are involved. It is intrinsic to evaluate them through the analytical expressions. However, as stated in \cite{EXP, EXP2, EXP3}, it is quite time-consuming for exponential function evaluation since numerous numerical operations are required to obtain convergence results. Unfortunately, even worse redundant calculations of those functions may exist when we fill the impedance matrix in MOM. To better understand this problem, we can roughly estimate the count of function calls to construct the impedance matrix as follows.

To obtain the quantitative analysis, we make the following assumptions. The structures are first discretized into $N$ planar triangles. $n$-point numerical integration at each edge of the source triangle to evaluate the inner surface integration is used and three lines overall for each triangle are required. $m$-point numerical integration on average for the outer surface integration is considered. Therefore, the total count of exponential function calls approximately equal to $3mnN^2$. It is easy to see that for a large number $N$ the count is extremely large because it depends on $N^2$. When acceleration techniques, like MLFMM/MLFMA \cite{MLFMP}, pFFT \cite{FastImp} and AIM \cite{AIM}, are used, the total count of function calls is approximately reduced to $3mnNC$, where $C$ is the averaged count of triangles in the near field region. We only consider the near-field interaction calculation. In those scenarios, in which the electrically large or multiscale structures are often involved, $N$ itself is an extremely large number. Therefore, we still make a large number of function calls to evaluate (\ref{EXPRR}) which is inevitably time-consuming. Further careful investigations can find that redundant function calls possibly exist in the direct evaluation method since (\ref{EXPRR}) are functions of distance between the field point or its project point and the source point.

In the following subsection, an adaptive interpolation method along with an efficient approach to retrieve sampling points is proposed to fast calculate (\ref{EXPRR}).

\subsection{Efficient Evaluations of $Exp$ Functions}
To efficiently evaluate values for the $Exp$ type functions with complex numbers, in the proposed method we should first determine sampling intervals. As shown in Fig. 1 (a) and (b), real and imaginary parts of $Exp(-jkr)$ and $Exp(-jkr)/r$ show different oscillating behaviors, which implies that special attention must be paid to selecting appropriate sampling points. $Exp(-jkr)$ is a periodic function without any singularity near the origin. Therefore, we can use a relatively uniform sampling interval and a linear interpolation method to obtain accurate interpolated values. For $Exp(-jkr)/r$, we start to sample it at a minimum distance $r_{min}$ near the origin since there is a singularity point where it shows large function variation. Large interpolation error can occur if the linear interpolation method is still used. In the proposed method, we select the Langrage polynomial method to interpolate $Exp(-jkr)/r$. We first uniformly sample the function using the interval $t$ defined as follows

\begin{equation}\label{INTERVAL}
{t} = \frac{{{\lambda _0}}}{{{{10}^4}\sqrt {{\varepsilon _r}{\mu _r}} }} \sim \frac{{{\lambda _0}}}{{{{10}^3}\sqrt {{\varepsilon _r}{\mu _r}} }},
\end{equation}
where $\lambda_0$,  $\epsilon_r$, $\mu_r$  denote the wavelength in free space, relative permittivity and relative permeability of homogenous medium, respectively. As stated in (\ref{INTERVAL}), to obtain acceptable accuracy of function values, we use $10^3 \sim 10^4$ sampling points per wavelength due to the error possibly accumulating in the integration phase when we evaluate double surface integrals. Since we only need to construct one-dimensional interpolation, the computational resource consumptions in terms of memory usage and CPU time can be ignored compared with those in the whole simulation even if we use such a dense sampling rate.
\begin{figure}
	\begin{minipage}[h]{0.48\linewidth}
		\centering
		\centerline{\includegraphics[width=1.95in]{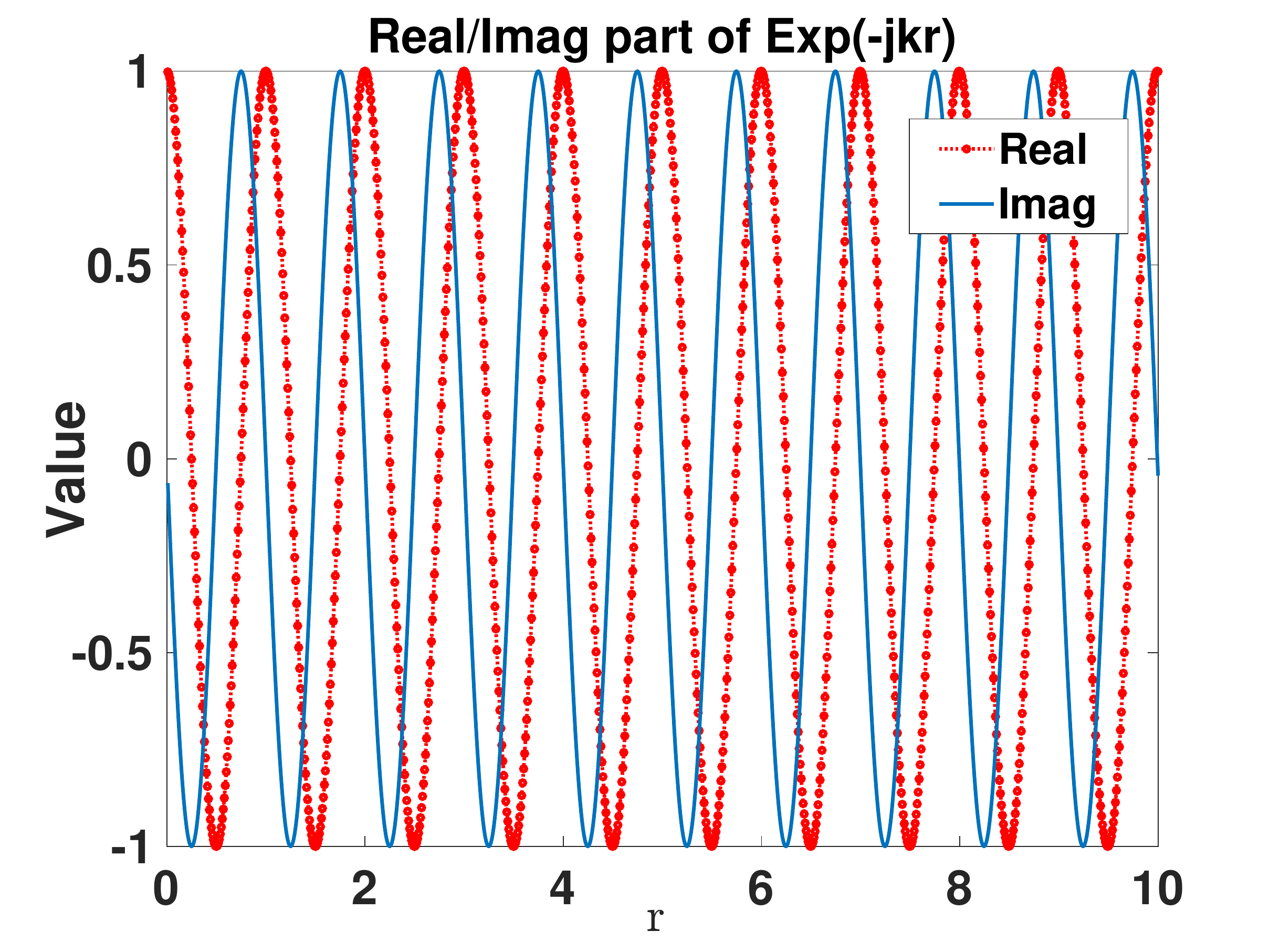}}		
		\centerline{(a)}
		\label{FIG1A}
	\end{minipage}
	\hfill
	\begin{minipage}[h]{0.48\linewidth}\label{FIG1B}
		\centerline{\includegraphics[width=1.95in]{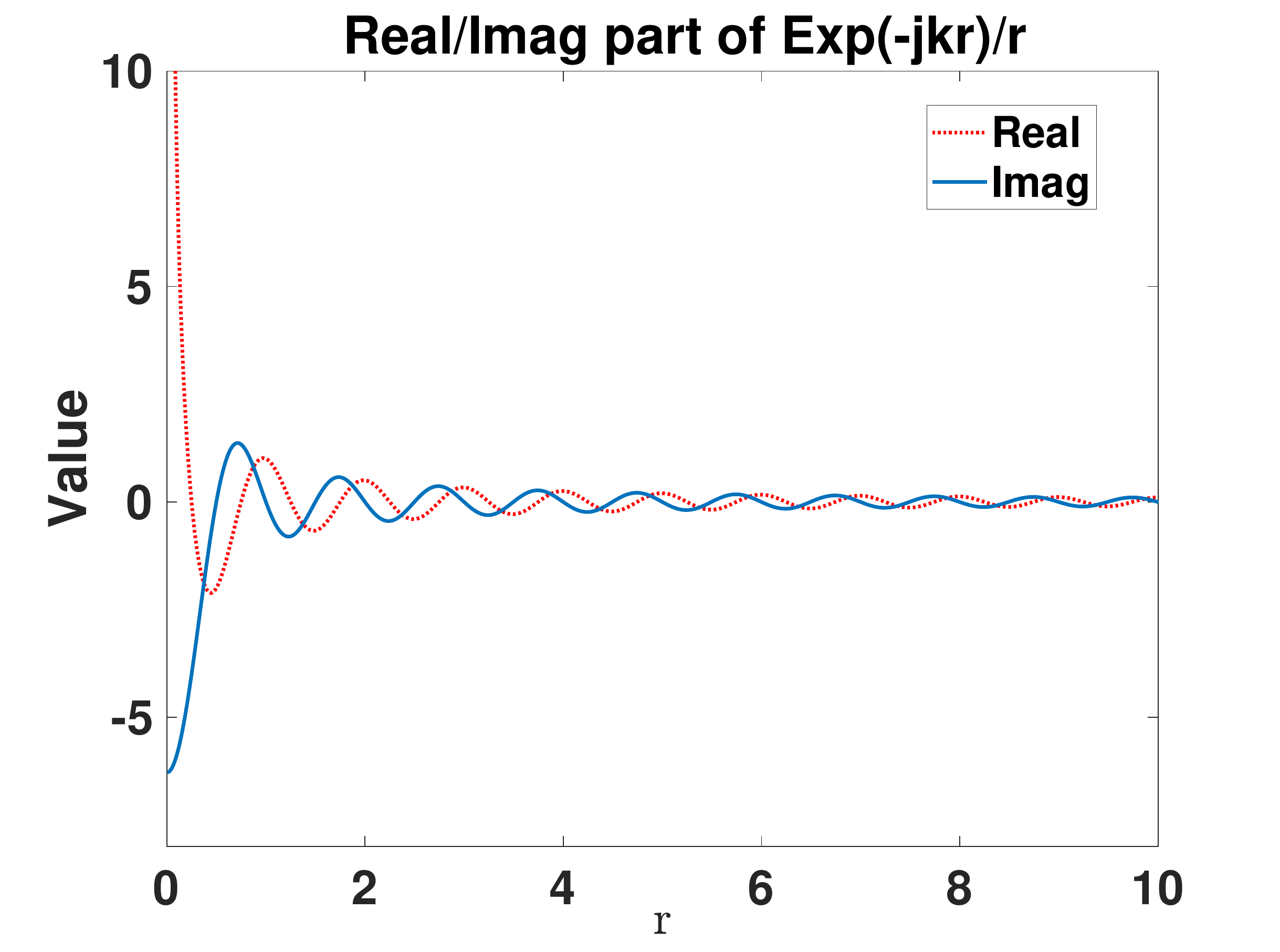}}
		\centerline{(b)}
	\end{minipage}
	\caption{Real and imaginary parts of $Exp(-jkr)$ and $Exp(-jkr)/r$.}
	\label{fig_6}
\end{figure}

However, we observed that quite large relative interpolation errors occur when the desired function values are small. We proposed the following method to keep the relative error uniformly distributed in the whole computational domain.

\begin{enumerate}
	\item Local refinement sampling is used near the zero function values. When sampling location is near zero function values, a local refinement with a smaller sampling interval is used.
	
	\item Zero function values are forced to be sampled. It is possible that the zero function value points are quite near those in $1)$. If it indeed happens, we only keep zero value sampled. This is because we want to keep $\Delta r_{min}$ reasonable large to make the hash array introduced in the next subsection with a reasonable dimension.
\end{enumerate}

Therefore, the proposed method can efficiently utilize the memory and construct an optimal interpolation tables. It is a banlance between the accuracy and computational cost. With the definition of the sampling points, several tables are constructed through their analytical expression without any redundant calculations. 

Various interpolation methods, like the linear interpolation method and the Lagrange interpolation method \cite{LINEAR}, are available to efficiently interpolate the continuous function values for (\ref{EXPRR}). They are used in the multilayer Green function evaluations \cite{AIM}\cite{MULTILAYERGREEN2}. The Lagrange interpolation method shows good performance in terms of accuracy and efficiency. Therefore, we use the Lagrange interpolation method in our study. However, it does not exclude other possibilities, such as the Bernstein, Bezier, and so on. Those polynomials should have similar performance in terms of accuracy.

We proposed to use the linear interpolation method to interpolate $Exp(-jkr)$ and the Lagrange polynomial method for $Exp(-jkr)/r$.  Therefore, the interpolated function values can be expressed as
\begin{equation}\label{INTERPOLATIONFunction}
f\left( {{r}} \right) = \sum\limits_{j = 0}^n {f\left( {{{{r}}_j}} \right){P_{n,j}}} ,
\end{equation}
where $r$ is the desired point location, $r_j$ is the sampling location, $n$ is the number of sampling points used in the interpolation, $P_{n,j}$ is the linear interpolation function 
\begin{equation}\label{INTERPOLATIONFunction}
P_{n,j}(r) = \frac{r-r_{j+1}}{r_j-r_{j+1}}~ \text{or}~ \frac{r-r_{j}}{r_{j+1}-r_{j}},
\end{equation}
or the Lagrange interpolation polynomial of degree $n$
\begin{equation}\label{INTERPOLATIONFunction}
P_{n,j}(r) = \prod \limits_{0 \le m \le k, m \ne j} \frac{{r - {r_m}}}{{{r_j} - {r_m}}}.
\end{equation}

\subsection{Efficient Searching Elements in the Sampling Table}
The dimension of sampling tables is usually quite large for electrically large or multiscale structures. If we naively search elements one-by-one in the sampling tables required by the interpolation methods, the performance will be severely deteriorated and even worse than that of original analytical method. We proposed the following efficient approach to solve the problem.

The hash-table method is designed to fast obtain the desired elements in a large container \cite{DATASTRUTURE}. We use the hash method to map $r$ to the index of sampling point in the sampling tables. The following many-to-one hash function is proposed
\begin{equation}\label{HashTable}
\left[\frac{{r - {r_{\min }}}}{{\Delta {r_{\min }}}}\right]
\end{equation}
where $r$ is the location of desired function values, $r_{min}$ is the minimum sampling location, $\Delta r_{min}$ is the minimum sampling interval, $\left[ \cdot \right] $ denotes the integer part of the corresponding float number. Once the desired $r$ is obtained we can fast locate the sampling point in the tables through (\ref{HashTable}). If $r \in [0,r_{min})$, we still use classical singularity cancellation method {\cite{GIBC}} to calculate the double surface integration. On the implementation aspect, one can simply construct a one dimensional array with the dimension of $[r_{max} - r_{min}]/(\Delta r_{min})$, where $r_{max}$ is the maximum sampling location. We first fill the array with the values of sampling locations using (\ref{HashTable}) and then other positions are recursively filled by its previous elements starting at the lower bound of the hash array. It can be used to fast retrieve the index of adaptively sampled function values stated in the previous subsection. Since the calculation of Green function at each integration point is independent of each other, the proposed method can be quite easy to be vectorized using vector CPU.

\section{Error Estimates}

\subsection{Error Estimates for Linear Interpolation}
Since in the proposed method we interpolate function values for (\ref{EXPRR}) in $[r_{min}, r_{max}]$, where they are infinitely smooth, the $n$-th derivative of these two functions are always available. The error bound for the polynomial interpolations for complex number can be easily derived according to the Mean Value theorem similar in \cite{LINEAR} as follows
\begin{equation}\label{ErrorLinearInt}
\mathop {\max }\limits_{r \in \left[ {a,b} \right]} \left| {f\left( {{r}} \right) - {f_n}\left( {{r}} \right)} \right| \le 
\mathop {\max }\limits_{r \in \left[ {a,b} \right]} \left| {{\omega _{n + 1}}\left( {{r}} \right)} \right|\frac{{\mathop {\max }\limits_{r \in \left[ {a,b} \right]} \left| {{f^{n + 1}}\left( {{\xi _r}} \right)} \right|}}{{\left( {n + 1} \right)!}},
\end{equation}
where
\begin{equation}
{\omega _{n + 1}}\left( {{r}} \right) = \prod\limits_{j = 0}^n {\left( {{{r}} - {{{r}}_j}} \right)} ,
\end{equation}
$\xi_r \in [a, b]$, $f(r)$ is the complex value function, $f_n(r)$ is the interpolated function value in $[a, b]$, $f^{n+1}$ is the $(n+1)^{th}$ derivative of $f$.

For the uniform sampling points, denoted as
\begin{equation}\label{UNIFORMSAMPLING}
{r_i} = a + \frac{i}{n}\left( {b - a} \right),i = 0,1,...,n,
\end{equation}
we can arrive at the analytical error bound with some mathematical manipulations based on (\ref{ErrorLinearInt}) as follows
\begin{equation}\label{ERRORUNIFORMSAMPLING}
\mathop {\max }\limits_{r \in \left[ {a,b} \right]} \left| {f\left( {{r}} \right) - {f_n}\left( {{r}} \right)} \right| \le {\left( {\frac{{b - a}}{n}} \right)^{n + 1}}\frac{\mathop {\max }\limits_{r \in \left[ {a,b} \right]} \left| {{f^{n + 1}}\left( {{\xi _r}} \right)} \right|}{{4\left( {n + 1} \right)}}.
\end{equation}

It is obvious that the error decreases as we increase the interpolation order or use smaller interpolation intervals. This is the reason why we locally refine the intervals near zero function values to obtain uniform relative error bound in the whole simulation domain. Therefore, to balance memory cost and accuracy, we use the uniform sampling far away from zero function values and non-uniform sampling near zero function values in the proposed method.

\subsection{Error Estimates for Lagrange Interpolation}
It is easy to derive the error bound for the Lagrange interpolation through the Generalized Rolle's theorem \cite{LINEAR} as follows
\begin{equation}\label{LAGRANGEINT}
\mathop {\max }\limits_{r \in \left[ {a,b} \right]} \left| {f\left( {{r}} \right) - {f_n}\left( {{r}} \right)} \right| \le \mathop \prod \limits_{j = 0}^n \left| {r - {r_j}}  \right| \frac{{\left| {{f^{n + 1}}\left( {{\xi _r}} \right)} \right|}} {{\left( {n + 1} \right)!}}.
\end{equation}
The general remarks for (\ref{LAGRANGEINT}) are the same as that of (\ref{ERRORUNIFORMSAMPLING}). The Lagrange interpolation shows better performance compared with the linear interpolation for highly oscillating kernel. However, the Lagrange method is more time-consuming compared with the linear interpolation. As stated before, we proposed the linear interpolation for $Exp(-jkr)$ and the Lagrange interpolation method for $Exp(-jkr)/r$ to balance the accuracy and efficiency.

\section{Numerical Experiments and Discussion}
\subsection{Accuracy of the proposed method}
We first validate the accuracy of the proposed method with an interpolation problem. The configurations are set up as follows, $r_{max} = 1$ m, corresponding to one wavelength at $300$ MHz, $r_{min} = 0.0001$ m, the sampling interval far away from zero function value is $t$, local refinement sampling interval is $t/2$, and local refinement range is $4t$. We interpolated $10^4$ function values uniformly distributed in $[0.0001~1]$. The $2$-point linear method and Lagrange polynominal method of order 3 are used in our simulations. The maximum relative error is defined as $\max \left( {{\left| {f\left( {{r}} \right) - {f_i}\left( {{r}} \right)} \right|}/{\left| {{f}\left( {{r}} \right)} \right|}} \right)$, where $f, f_i$ are the analytical function value and the interpolated function value, respectively.

\begin{figure}
	\begin{minipage}[h]{0.48\linewidth}\label{FIG2A}
		\centering
		\centerline{\includegraphics[width=1.95in]{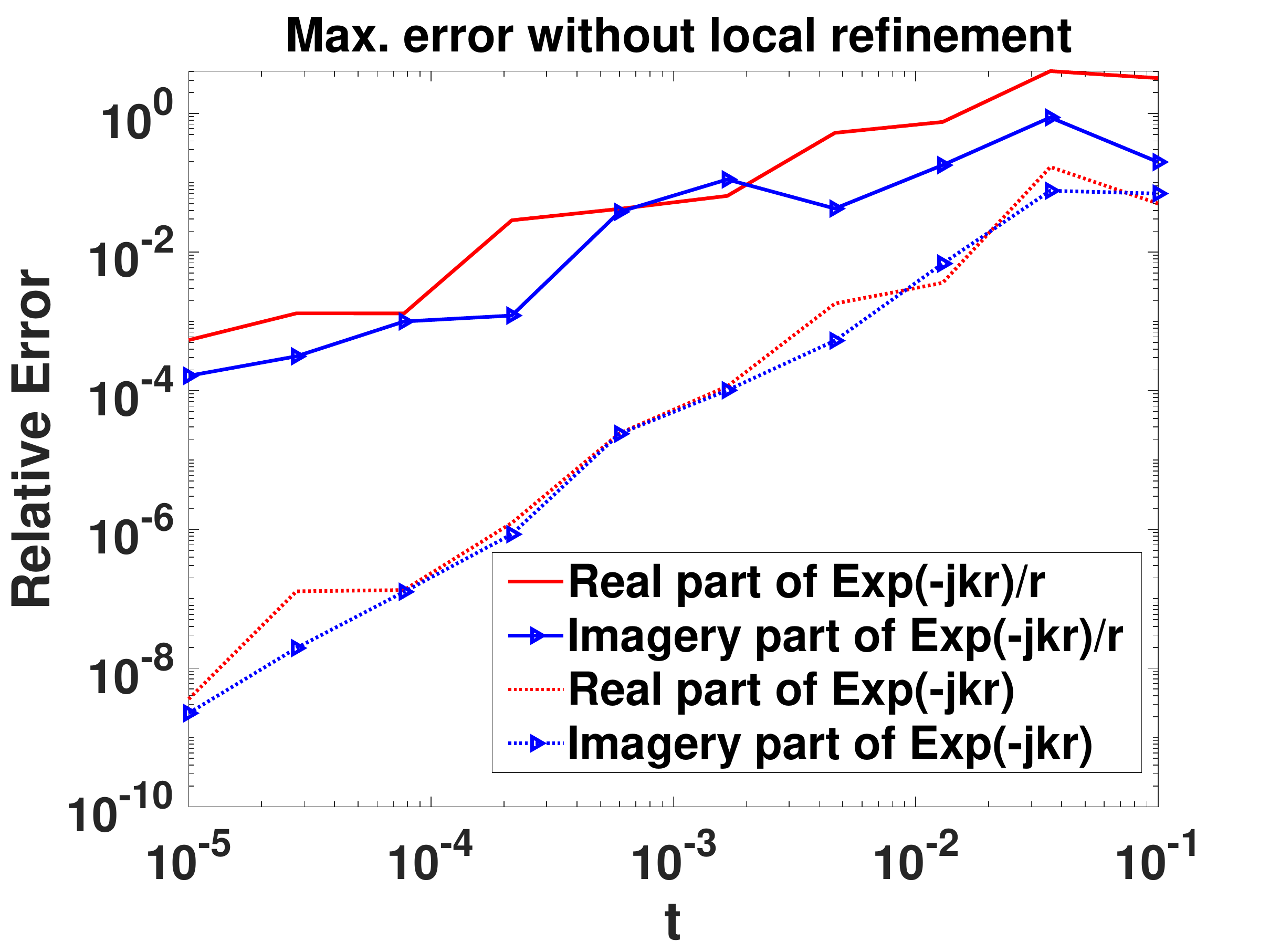}}		
		\centerline{(a)}
	\end{minipage}
	\hfill
	\begin{minipage}[h]{0.48\linewidth}\label{FIG2B}
		\centerline{\includegraphics[width=1.95in]{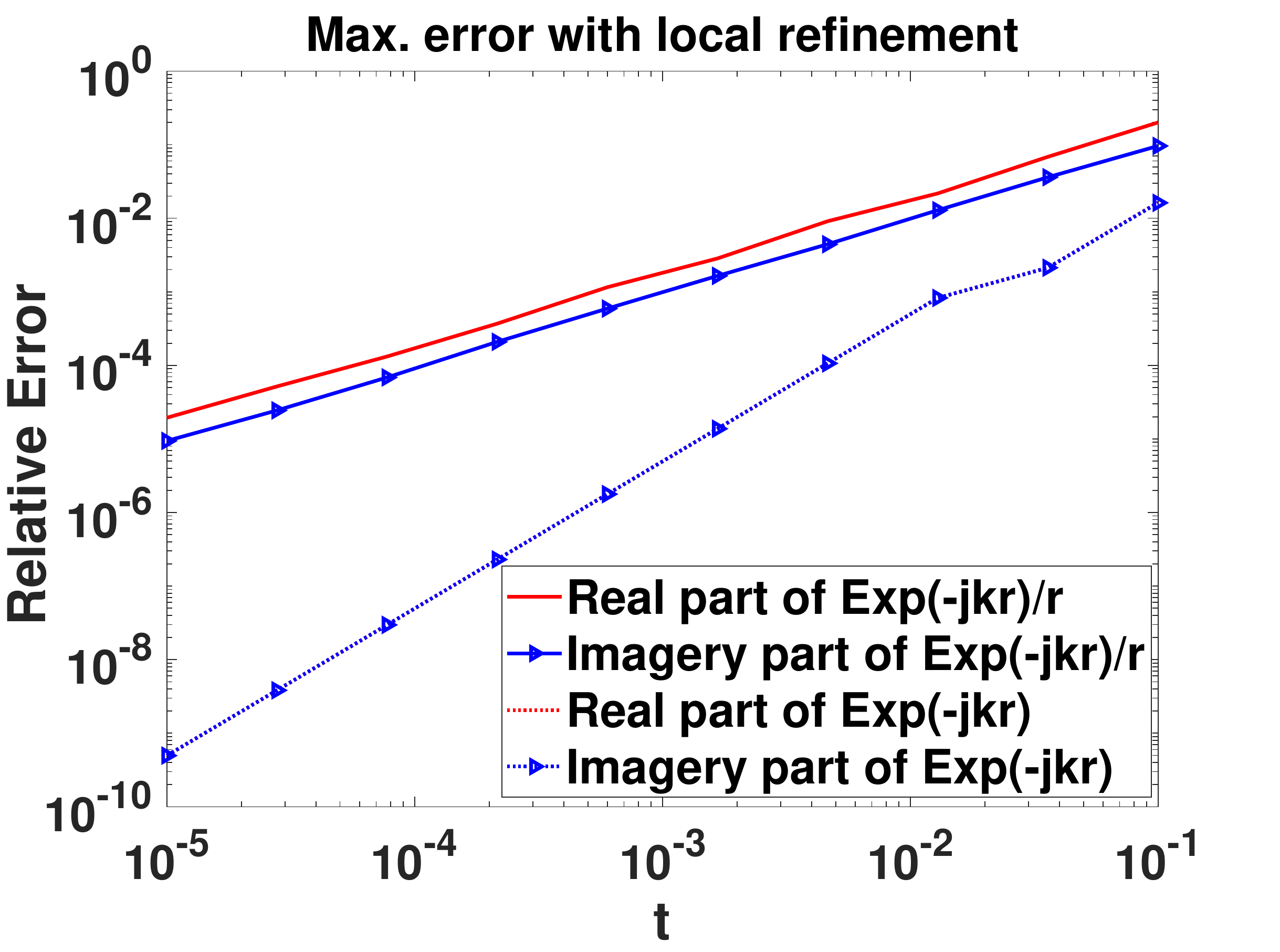}}
		\centerline{(b)}
	\end{minipage}
	\hfill
	\begin{minipage}[h]{0.48\linewidth}\label{FIG2C}
		\centerline{\includegraphics[width=1.95in]{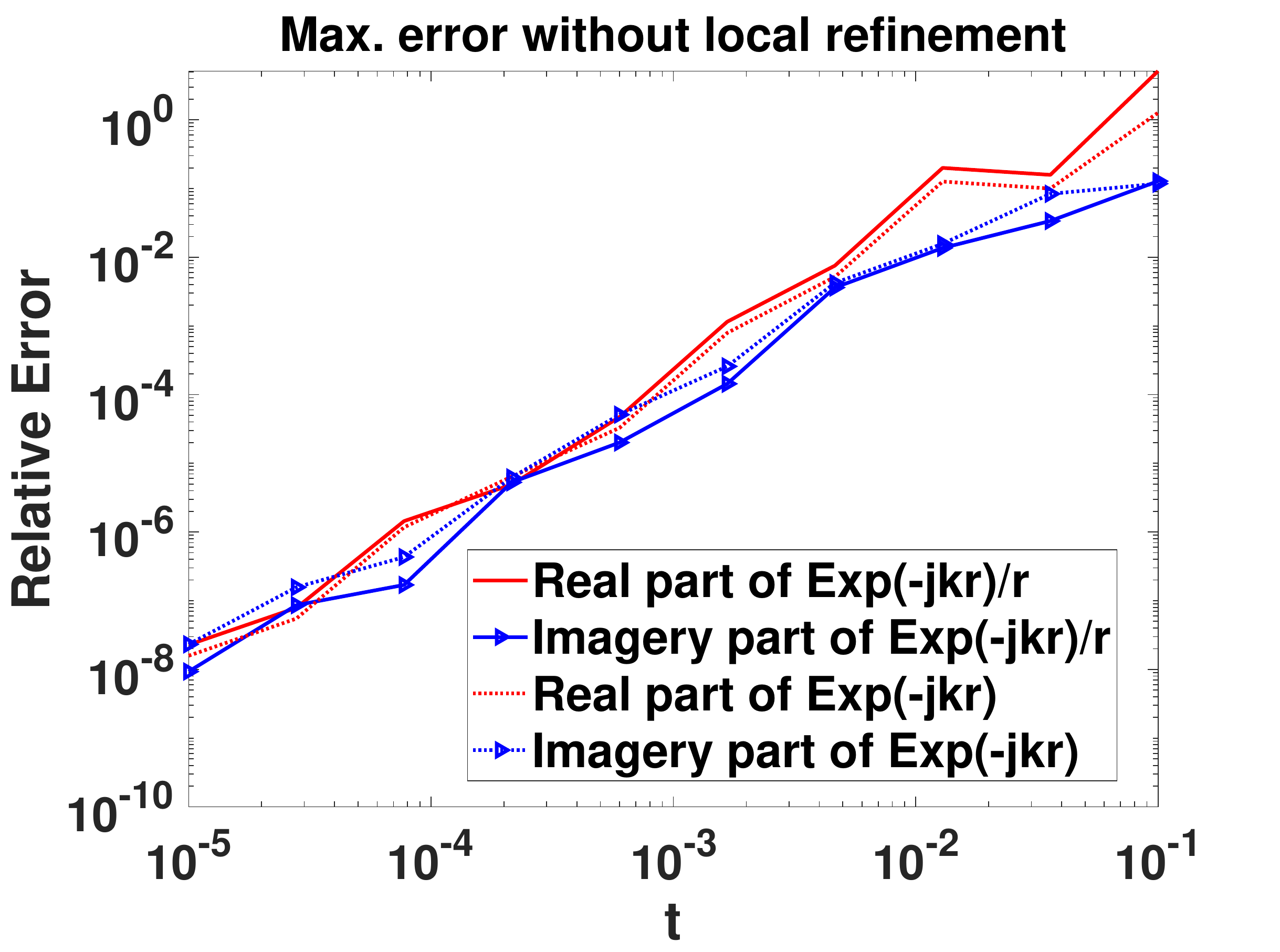}}
		\centerline{(c)}
	\end{minipage}
	\hfill
	\begin{minipage}[h]{0.48\linewidth}\label{FIG2D}
	\centerline{\includegraphics[width=1.95in]{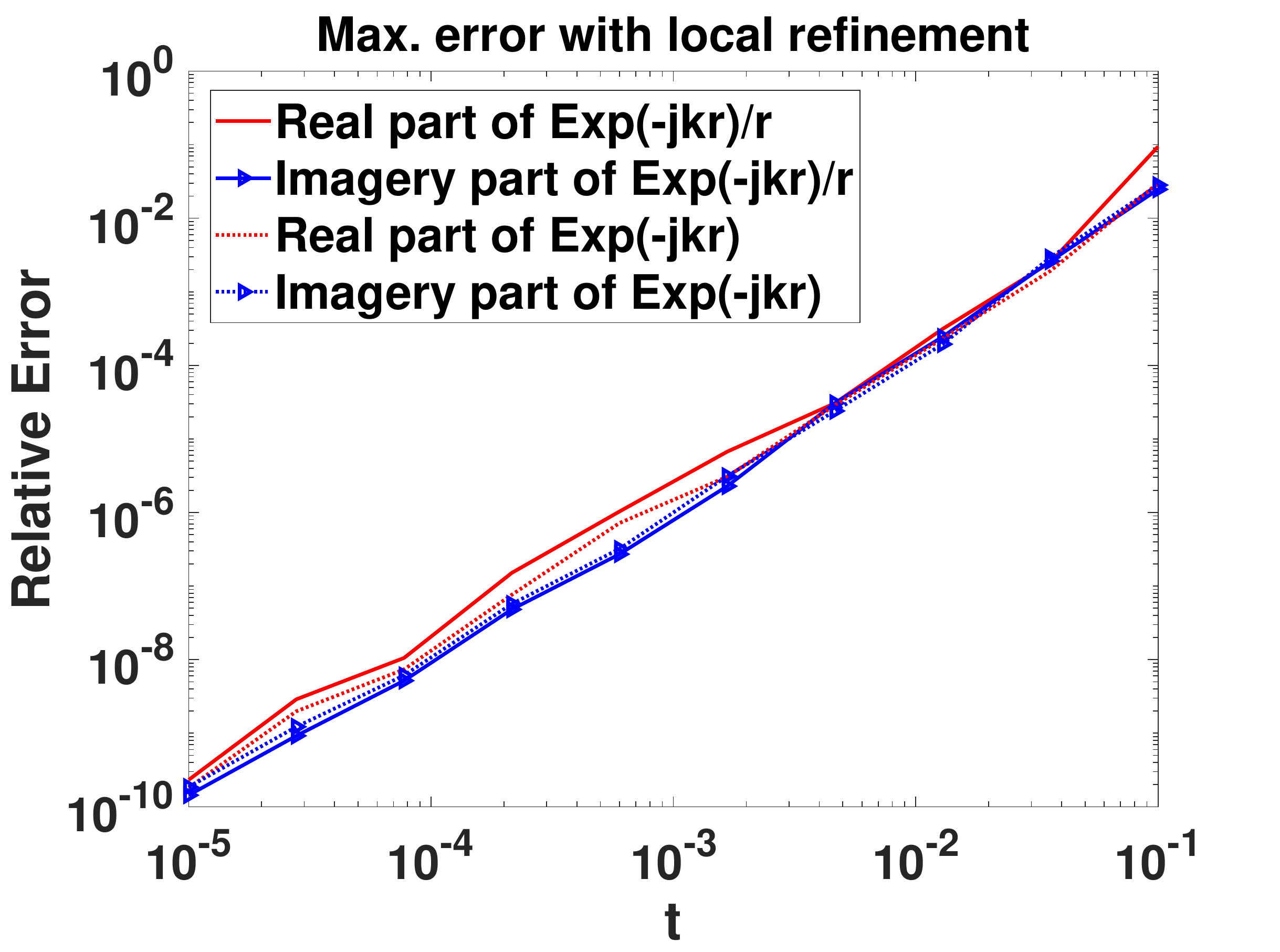}}
	\centerline{(d)}
	\end{minipage}
	\caption{(a) and (b) denote the maximum relative error of real and imaginary part of $Exp(-jkr)$ and $Exp(-jkr)/r$ obtained from linear interpolation method, (c) and (d) from Lagrange polynomial method. }
	\label{fig_6}
\end{figure}

It is easy to see that in Fig. 2 (a) and (b), (c) and (d), the proposed local refinement method can approximately improve the accuracy up to two orders for both the real and imaginary parts of (\ref{EXPRR}). Compared Fig. 2 (a) and (b) with (c) and (d), the linear method shows similar performance in terms of the maximum relative error for $Exp(-jkr)$ as that of the Lagrange polynomial method. Since the linear method has higher efficiency than that of Lagrange, we use the linear method to interpolate $Exp(-jkr)$ in the proposed method.

As shown in Fig. 2, the Lagrange polynomial method is significantly more accurate than the linear method for $Exp(-jkr)/r$. Therefore, we use the Lagrange method for the Green function to guarantee the accuracy of the proposed method.

\begin{figure}
	\begin{minipage}[h]{0.48\linewidth}\label{FIG2A}
		\centering
		\centerline{\includegraphics[width=1.85in]{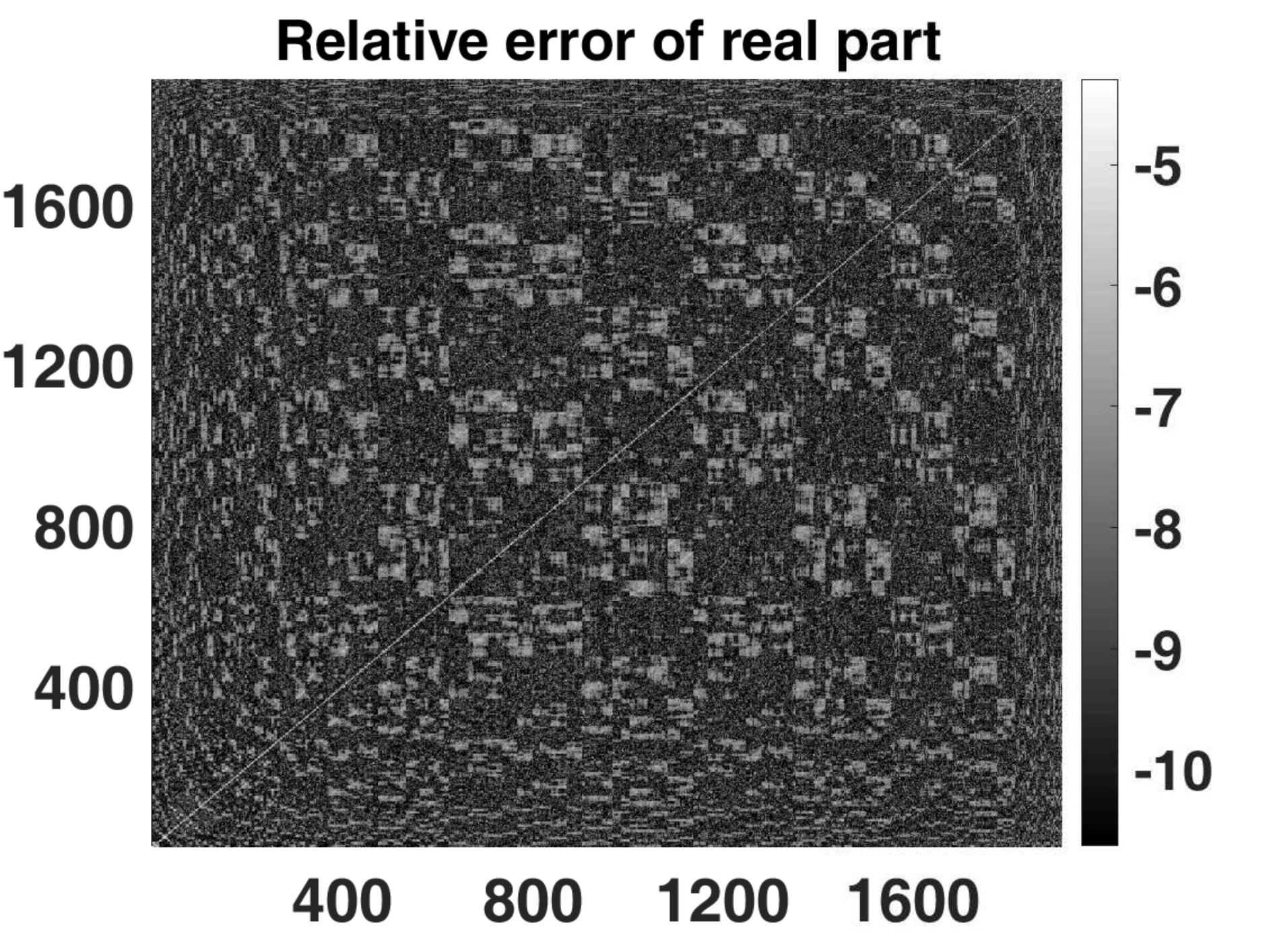}}		
		\centerline{(a)}
	\end{minipage}
	\hfill
	\begin{minipage}[h]{0.48\linewidth}\label{FIG2B}
		\centerline{\includegraphics[width=1.85in]{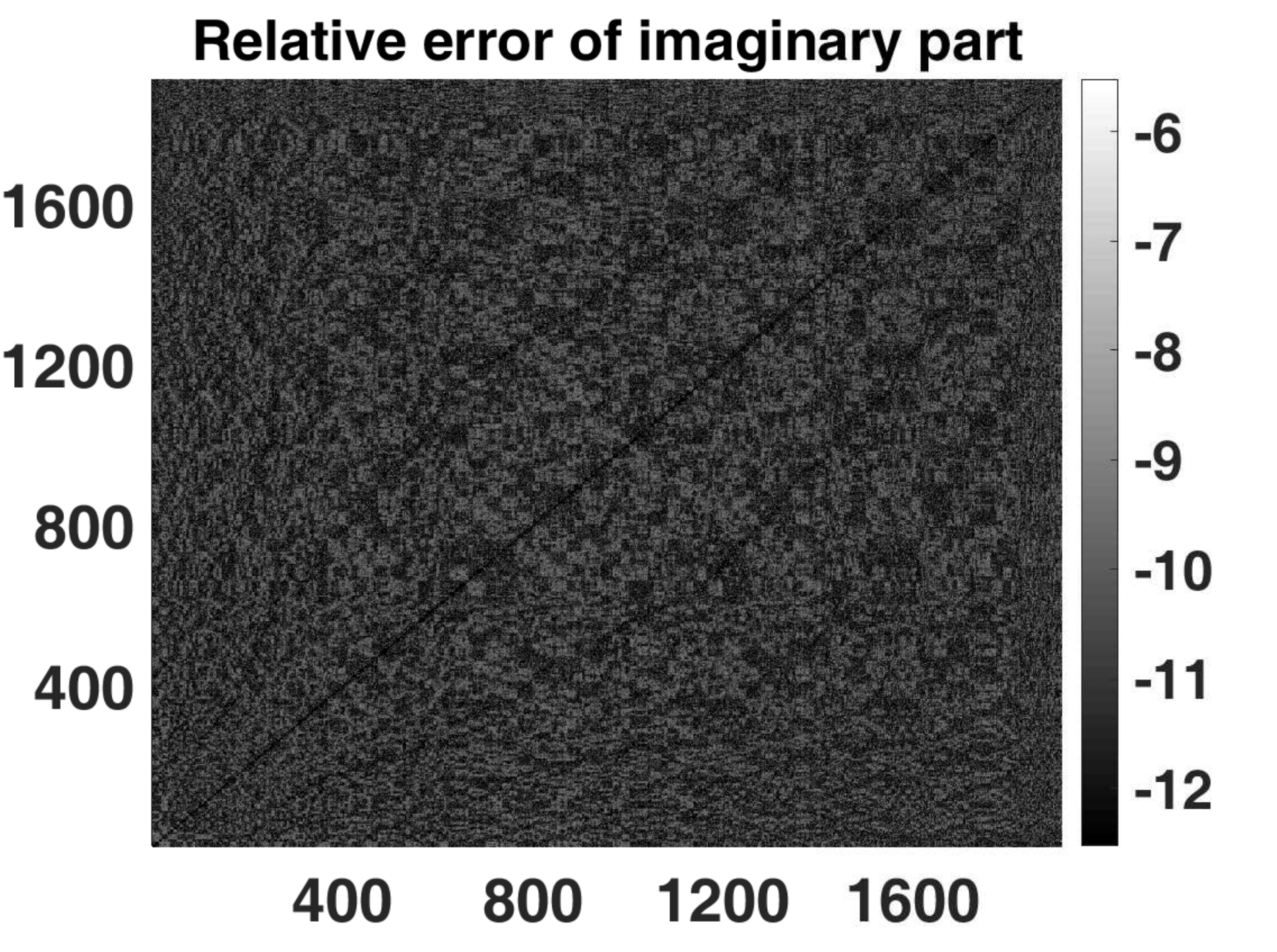}}
		\centerline{(b)}
	\end{minipage}
	\hfill
	\begin{minipage}[h]{0.48\linewidth}\label{FIG2C}
		\centerline{\includegraphics[width=1.85in]{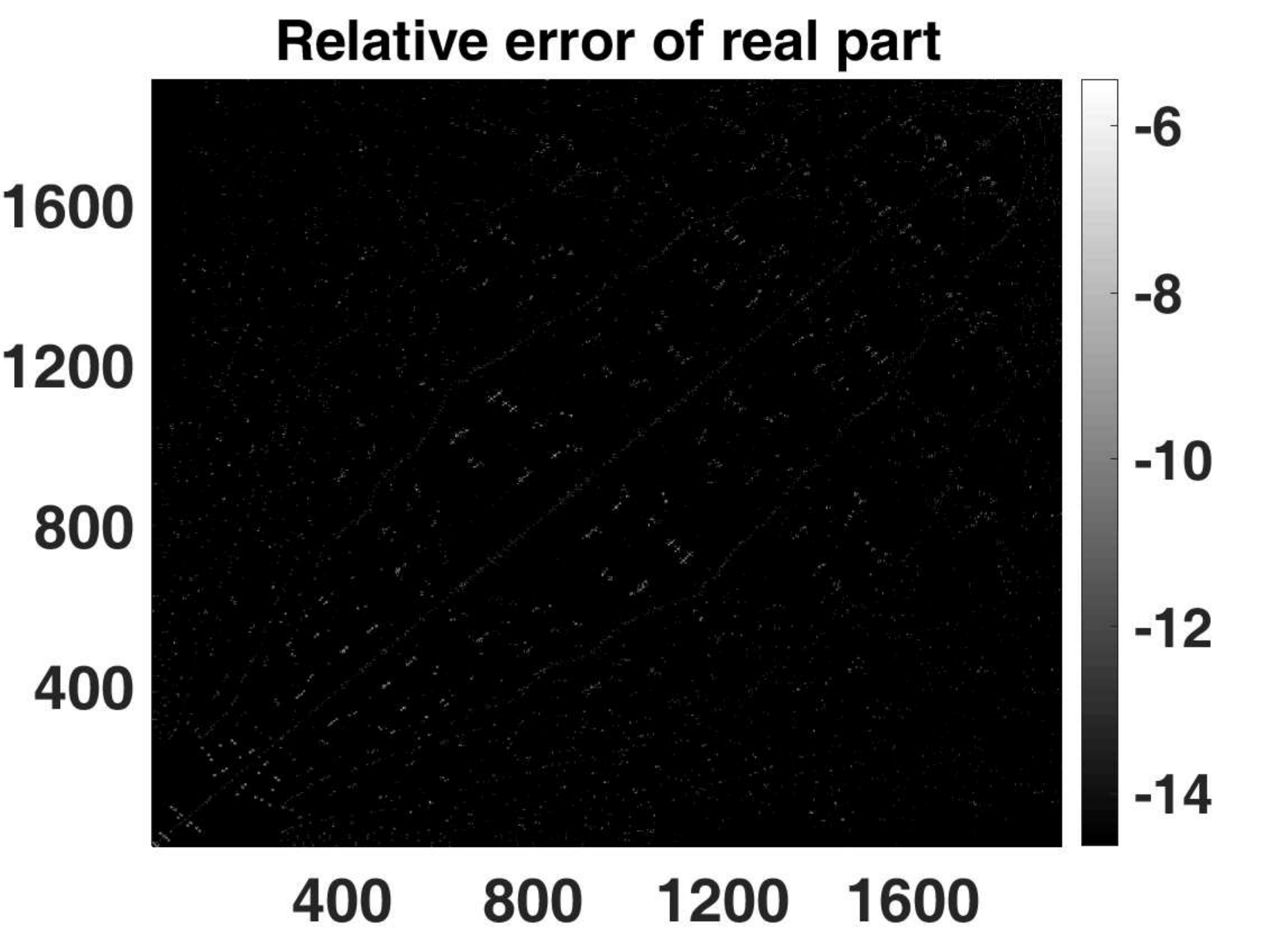}}
		\centerline{(c)}
	\end{minipage}
	\hfill
	\begin{minipage}[h]{0.48\linewidth}\label{FIG2D}
		\centerline{\includegraphics[width=1.85in]{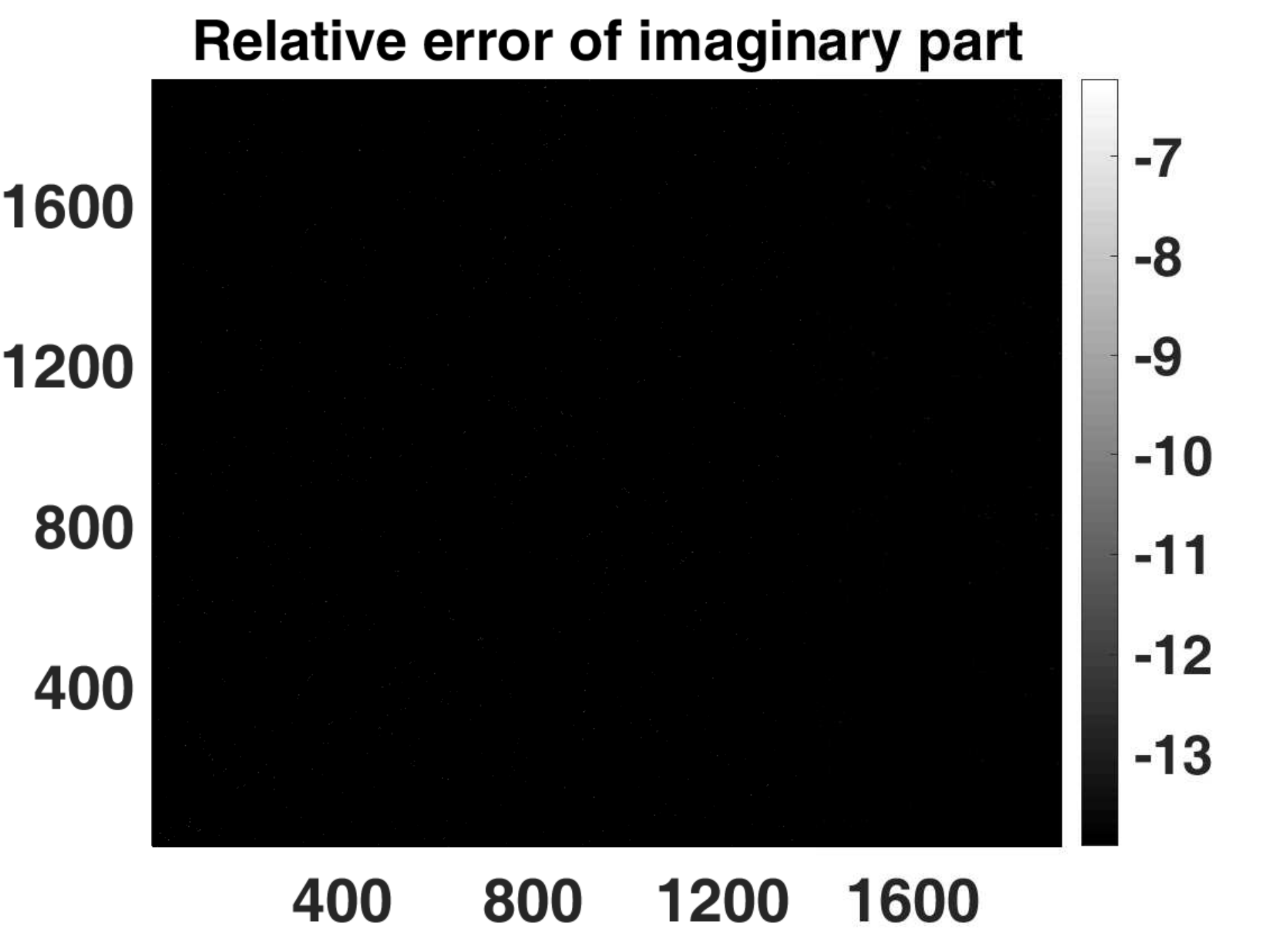}}
		\centerline{(d)}
	\end{minipage}
	\caption{The relative error of real and imaginary part of elements in impedance matrix through EFIE and MFIE obtained from the proposed method. Note: the relative errorrs are in log scale.}
	\label{fig_6}
\end{figure}

\subsection{Performance of the Proposed Method Incorporated in MOM based on Various Formulations}
A sphere with the radii of 0.5 wavelength and perfectly electrical conductor (PEC) is considered for the EFIE and the MFIE.  A $x$-polarized plane wave incidents along $z$ axis. All numerical experiments are carried out on a workstation with an intel i7-7700 3.6 GHz CPU and 32 G memory and our in-house codes are running with a thread without exploring any the parallel computation to make fair comparisons.

In Fig. 3, we plot the relative error of all elements in the impedance matrix obtained from the EFIE involved $L$ operator and the MFIE involved $K$ operator \cite{MOMBOOK}. Other formulations in essence are various combinations of the EFIE and MFIE, therefore, the performance of the proposed method in other formulations are similar. As shown in Fig. 3 (a) and (b), (c) and (d), the relative error of both the real and imaginary parts for the EFIE reaches $10^{-5}$ and $10^{-6}$. For the MFIE, the relative error is further decreased to $10^{-6}$ and $10^{-7}$. That implies the proposed method can provide an excellent level of elementwise accuracy. The $L_2$ error of RCS for both EFIE and MFIE can reach $10^{-6}$. Those results show excellent performance in terms of accuracy for the proposed method.

\begin{table}[!t]
	\begin{center}
		\caption{COMPARISON OF TIME COST OF DIFFERENT FORMULATIONS BETWEEN THE ANALYTICAL METHOD AND THE PROPOSED METHOD}
		\begin{tabular}{c c c c c }		
			\toprule
			\midrule
			\multirow{2}[4]{*}{Metric} & \multicolumn{2}{c}{EFIE }  & \multicolumn{2}{c}{MFIE }\\
			\cmidrule(rl){2-5}
			& Analytical                               & Proposed & Analytical & Proposed  \\
			\cmidrule(r){1-2}\cmidrule(l){2-5}
			\multicolumn{1}{l}{Counts of unknowns}     & 1914   & 1914  & 1914   & 1914  \\
			\multicolumn{1}{l}{Matrix filling time [s]}& 85     & 61    & 82   & 62  \\
			\multicolumn{1}{l}{Overall time [s]}       & 87     & 63    & 84   & 64  \\
			\midrule
			\bottomrule
		\end{tabular}
	\end{center}
\end{table}

Table I shows the time cost for the EFIE and the MFIE incorporated with the analytical method and the proposed method. For the EFIE, matrix filling in the analytical method takes 85s. However, when incorporated with the proposed method, matrix filling only takes 61s with $L_2$ error of RCS  $10^{-6}$. It shows 27\% and 24 \% efficiency improvement for the EFIE and the MFIE formulations, respectively, without sacrificing accuracy. 

Another practical inductor \cite{INDUCTORS} used in the ICs has been modelled through the MOM with and without  the proposed method. It was discretized with 728 triangles and 1092 edges. We have extracted the scattering parameters from 1 GHz to 30 GHz. The relative error at each frequency point reaches the order of $10^{-3}$. The simulation without the proposed method takes  405 seconds. However, when our MOM codes are incorporated with the proposed method, it only takes 315 seconds to complete the simulation. It shows the proposed method can save 22\% computational time without sacrificing the accuracy. Therefore, the proposed method is also applicable to the practical structures and show excellent performance in terms of  efficiency and accuracy.

\section{Conclusions}
In this paper, we proposed an efficient interpolation approach to accelerate matrix filling in MOM. In the proposed method, several interpolation tables are constructed with carefully consideration of local refinement near the zero function values. The linear method and Lagrange polynomial are used to interpolate $Exp(-jkr)$ and $Exp(-jkr)/r$, respectively. An efficient approach to obtain sampling points are also presented. Compared with the analytical method, it shows over 20\% efficiency improvement without sacrificing accuracy. In all, the proposed method is simple, robust and efficient to accelerate matrix filling for MOM and applicable for almost all SIE formulations and it is quite easy to be integrated into the available MOM codes with only a few modifications.



%

\end{document}